\documentclass[12pt]{article}
\usepackage{amssymb}
\usepackage{amsbsy}
\begin{document}
\centerline{\Large Complex-valued valuations on $L^{p}$ spaces} \vskip 4pt

\begin{center}
{~~Lijuan~Liu~~
\vskip 15pt

School of Mathematics and Computational Science,
Hunan University of Science and Technology, Xiangtan, 411201, P.R.China. \\E-mail: lijuanliu@hnust.edu.cn
}

\end{center} \vskip 5pt

\begin{center}
\begin{minipage}{13cm}
{\bf Abstract:}
All continuous translation invariant complex-valued valuations on Lebesgue measurable functions are completely classified.
And all continuous rotation invariant complex-valued valuations on spherical Lebesgue measurable functions are also completely classified.

\vskip 15pt {\bf MR(2010)  Subject Classification}~~  52B45, 52A20
\end{minipage}
\end{center}

 \vskip 5pt
 \vskip 10pt \centerline{\bf \S1.~~Introduction}

\vskip 15pt
\noindent A function $z$ defined on a lattice $(\mathcal{L},\vee,\wedge)$ and taking values in an abelian semigroup is called a \textit{valuation} if
$$z(f\vee g)+z(f\wedge g)=z(f)+z(g)\eqno(1.1)$$
for all $f,g\in \mathcal{L}$. A function $z$ defined on some subset $\mathcal{L}_{0}$ of $\mathcal{L}$ is called a valuation on $\mathcal{L}_{0}$ if
(1.1) holds whenever $f, g, f\vee g, f\wedge g\in \mathcal{S}$. For $\mathcal{L}_{0}$ the set of convex bodies, $\mathcal{K}^{n}$, in
$\mathbb{R}^{n}$ with $\vee$ denoting union and $\wedge$ intersection. Valuations on convex bodies are a classical concept going back to Dehn's solution in 1900 of Hilbert's Third Problem.
Probably the most famous result on valuations is Hadwiger's classification theorem of continuous rigid motion invariant valuations. For the more recent contributions on valuations convex bodies
can see [1,2,4,6,12-24,26-33,38,39,43-45].

Valuations on convex bodies can be
considered as valuations on suitable function spaces. Recently, valuations on functions has been rapidly growing (see [3,5,7-11,25,34-37,40-42,46-49]). For a space of real-valued functions, the operations $\vee$ and $\wedge$ are defined as pointwise maximum and minimum, respectively.
A complete classification of valuations intertwining with the SL(n) on Sobolev spaces [34,36,40] and $L^{p}$ spaces [25,35,37,42,46,47,49] was established, respectively. Valuations on convex functions [3,7,10,11], log-concave functions [41], quasi-concave functions [8,9] and functions of Bounded variations [48] were introduced and classified.

Recently, Wang and the author [49] showed that the Fourier transform is the only valuation which is a continuous, positively GL(n) covariant and logarithmic translation covariant complex-valued valuation on integral functions. This motivates the study of complex-valued valuations on functions.

Let $\mathcal{L}$ be a lattice of complex-valued functions.
For $f\in \mathcal{L}$, let $\Re f$ and $\Im f$
denote the real and imaginary parts of $f$, respectively.
The pointwise maximum of $f$ and $g$, $f\vee g$ and the pointwise minimum $f\wedge g$ are defined by
$$f\vee g=\Re f\vee\Re g+i(\Im f\vee\Im g),\eqno(1.2)$$
and
$$f\wedge g=\Re f\wedge\Re g+i(\Im f\wedge\Im g).\eqno(1.3)$$ If $f,g$ are real-valued functions, then (1.2) and (1.3) coincide with the real cases.
A function $\Phi: \mathcal{L}\rightarrow \mathbb{C}$ is called a \textit{valuation} if
$$\Phi(f\vee g)+\Phi(f\wedge g)=\Phi(f)+\Phi(g)$$
for all $f,g\in \mathcal{L}$ and $\Phi(0)=0$ if $0\in \mathcal{L}$.
It is called \textit{continuous} if
$$\Phi(f_{i})\rightarrow \Phi(f), ~\textrm{as}~~f_{i}\rightarrow f ~\textrm{in}~\mathcal{L}.$$
It is called \textit{translation invariant} if
$$\Phi(f(\cdot-t))=\Phi(f)$$
for every $t\in \mathbb{R}^{n}$. It is called \textit{rotation invariant} if
$$\Phi(f\circ \theta^{-1})=\Phi(f)$$
for every $\theta\in \textrm{O}(n)$, where $\theta^{-1}$ denotes the inverse of $\theta$.

Let $p\geq 1$. If $(X,\mathfrak{F},\mu)$ is a measure space, then the space $L^{p}$-space, $L^{p}(\mathbb{C},\mu)$ is the collection of $\mu$-measurable complex-valued
functions $f: X\rightarrow \mathbb{C}$ that satisfies
$$\int_{X}|f|^{p}d\mu<\infty.$$
A measure space $(X,\mathfrak{F},\mu)$
is called \textit{non-atomic} if for every $E\in\mathfrak{F}$ with $\mu(E)>0$, there exists $F\in\mathfrak{F}$ with $F\subseteq E$ and $0<\mu(F)<\mu(E)$. Let $\chi_{E}$ denote the characteristic function of the measurable set $E$, i.e.
$$\chi_{E}(x)=\left\{
                \begin{array}{ll}
                  1, & \hbox{$x \in E$;} \\
                  0, & \hbox{$x\notin E$.}
                \end{array}
              \right.
$$
\vskip 5pt \noindent {\bf Theorem 1.1} Let $(X,\mathfrak{F},\mu)$ be a non-atomic measure space and let $\Phi: L^{p}(\mathbb{C},\mu)\rightarrow \mathbb{C}$ be a continuous valuation. If
there exist continuous functions $h_{k}:\mathbb{R}\rightarrow \mathbb{R}$ with $h_{k}(0)=0~(k=1,2,3,4)$ such that
$\Phi(c\chi_{E})=(h_{1}(\Re c)+h_{3}(\Im c))\mu(E)+i(h_{2}(\Re c)+h_{4}(\Im c))\mu(E)$
for all $c\in \mathbb{C}$ and all $E\in\mathfrak{F}$ with $\mu(E)<\infty$, then there exist constants $\gamma_{k},\delta_{k}\geq 0$
such that $|h_{k}(a)|\leq \gamma_{k}|a|^{p}+\delta_{k}$ for $a\in \mathbb{R}$, and
$$\Phi(f)=\int_{X}(h_{1}\circ\Re f+h_{3}\circ\Im f)d\mu+i\int_{X}(h_{2}\circ\Re f+h_{4}\circ\Im f)d\mu$$
for all $f\in L^{p}(\mathbb{C},\mu)$. In addition, $\delta_{k}=0~(k=1,2,3,4)$ if $\mu(X)=\infty$.

Let $L^{p}(\mathbb{C},\mathbb{R}^{n})$ denote the $L^{p}$ space of Lebesgue measurable complex-valued functions on $\mathbb{R}^{n}$.
\vskip 5pt \noindent {\bf Theorem 1.2} A function $\Phi: L^{p}(\mathbb{C},\mathbb{R}^{n})\rightarrow \mathbb{C}$ is a continuous translation invariant valuation if and only if
there exist continuous functions $h_{k}:\mathbb{R}\rightarrow \mathbb{R}$ with the property that there exist constants $\gamma_{k}\geq 0$ such that
$|h_{k}(a)|\leq \gamma_{k}|a|^{p}$ for all $a\in\mathbb{R}~(k=1,2,3,4)$, and
$$\Phi(f)=\int_{\mathbb{R}^{n}}(h_{1}\circ\Re f+h_{3}\circ\Im f)(x)dx+i\int_{\mathbb{R}^{n}}(h_{2}\circ\Re f+h_{4}\circ\Im f)(x)dx\eqno(1.4)$$
for all $f\in L^{p}(\mathbb{C},\mathbb{R}^{n})$.

Let $S^{n-1}$ be the unit sphere in $\mathbb{R}^{n}$ and let $L^{p}(\mathbb{C},S^{n-1})$ denote the $L^{p}$ space of spherical Lebesgue measurable complex-valued functions on $S^{n-1}$.

\vskip 5pt \noindent {\bf Theorem 1.3} A function $\Phi: L^{p}(\mathbb{C},S^{n-1})\rightarrow \mathbb{C}$ is a
continuous rotation invariant valuation if and only if there exist continuous functions $h_{k}:\mathbb{R}\rightarrow \mathbb{R}$ with the properties that $h_{k}(0)=0$ and there exist constants $\gamma_{k}, \delta_{k}\geq 0$ such that
$|h_{k}(a)|\leq \gamma_{k}|a|^{p}+\delta_{k}$ for all $a\in\mathbb{R}~(k=1,2,3,4)$, and
$$\Phi(f)=\int_{S^{n-1}}(h_{1}\circ\Re f+h_{3}\circ\Im f)(u)du+i\int_{S^{n-1}}(h_{2}\circ\Re f+h_{4}\circ\Im f)(u)du\eqno(1.5)$$
for all $f\in L^{p}(\mathbb{C},S^{n-1})$.

 \vskip 5pt
 \vskip 10pt \centerline{\bf \S2.~~ Notation and preliminary results}

\vskip 15pt \noindent  We collect some properties of complex-valued functions. If $f$ is a complex-valued function on $\mathbb{R}^{n}$, then $$f(x)=\Re f+i\Im f,$$ where $\Re f$ and $\Im f$ denote the real part and imaginary part of $f$, respectively.
The absolute value of $f$ which is also called modulus is defined by $$|f|=\sqrt{(\Re f)^{2}+(\Im f)^{2}}.$$
Let $p\geq 1$. For a measure space $(X,\mathfrak{F},\mu)$, define $L^{p}(\mathbb{C},\mu)$ as the space of $\mu$-measurable complex-valued
functions $f: X\rightarrow \mathbb{C}$ that satisfies
$$\|f\|_{p}=(\int_{X}|f|^{p}d\mu)^{\frac{1}{p}}<\infty.$$
Let $f, g\in L^{p}(\mathbb{C},\mu)$, then $f\vee g, f\wedge g\in L^{p}(\mathbb{C},\mu)$. The functional $\parallel \cdot\parallel: L^{p}(\mathbb{C},\mu)\rightarrow \mathbb{R}$ is a semi-norm. If functions in $L^{p}(\mathbb{C},\mu)$ that are equal almost everywhere with respect to $\mu~ (\textrm{a.e.}~[\mu])$ are identified, then $\parallel \cdot\parallel: L^{p}(\mathbb{C},\mu)\rightarrow \mathbb{R}$ becomes a norm. Obviously, $L^{p}(\mathbb{C},\mu)$ is a lattice of complex-valued functions.
Let $L^{p}(\mathbb{R},\mu)$ denote the subset of $L^{p}(\mathbb{C},\mu)$ where the functions take real values.
For $f_{i}, f\in L^{p}(\mathbb{C},\mu)$, if $\|f_{i}-f\|_{p}\rightarrow0$, then $f_{i}\rightarrow f$ in $L^{p}(\mathbb{C},\mu)$. Moreover,
$$f_{i}\rightarrow f \in L^{p}(\mathbb{C},\mu) \Leftrightarrow \Re f_{i}\rightarrow \Re f, \Im f_{i}\rightarrow \Im f \in L^{p}(\mathbb{R},\mu).$$

The following characterizations of real-valued valuations on $L_{p}$ spaces which were established by Tsang [46] will play key role in our proof.
Let $L^{p}(\mathbb{R}, \mathbb{R}^{n})$ and $L^{p}(\mathbb{R}, S^{n-1})$ denote
the $L_{p}$ space of Lebesgue measurable real-valued functions on $\mathbb{R}^{n}$ and the
$L_{p}$ space of spherical Lebesgue measurable real-valued functions on $S^{n-1}$, respectively.

\vskip 5pt
\noindent {\bf Theorem 2.1} ([46]) A function $\Phi: L^{p}(\mathbb{R},\mathbb{R}^{n})\rightarrow \mathbb{R}$ is a continuous translation invariant valuation if and only if
there exists a continuous function $h:\mathbb{R}\rightarrow \mathbb{R}$ with the property that there exists a constant $\gamma\geq 0$ such that
$|h(a)|\leq \gamma|a|^{p}$ for all $a\in\mathbb{R}$, and
$$\Phi(f)=\int_{\mathbb{R}^{n}}(h\circ f)(x)dx$$
for all $f\in L^{p}(\mathbb{R},\mathbb{R}^{n})$.

\vskip 5pt \noindent {\bf Theorem 2.2} ([46]) A function $\Phi: L^{p}(\mathbb{R},S^{n-1})\rightarrow \mathbb{R}$ is a
continuous rotation invariant valuation if and only if there exists a continuous function $h:\mathbb{R}\rightarrow \mathbb{R}$ with the properties that $h(0)=0$ and there exist constants $\gamma, \delta\geq 0$ such that
$|h(a)|\leq \gamma|a|^{p}+\delta$ for all $a\in\mathbb{R}$, and
$$\Phi(f)=\int_{S^{n-1}}(h\circ f)(u)du$$
for all $f\in L^{p}(\mathbb{R},S^{n-1})$.

 \vskip 10pt \centerline{\bf \S3.~~Main results}
\vskip 10pt
\noindent {\bf Lemma 3.1} Let $h_{k}: \mathbb{R}\rightarrow \mathbb{R}$ be continuous functions with the properties that $h_{k}(0)=0$ and
there exist $\gamma_{k}, \delta_{k}\geq 0$ such that $|h_{k}(a)|\leq \gamma_{k}|a|^{p}+\delta_{k}$ for all $a\in\mathbb{R} (k=1,\ldots,4)$. If the function $\Phi: L^{p}(\mathbb{C},\mu)\rightarrow \mathbb{C}$ is defined by
$$\Phi(f)=\int_{X}(h_{1}\circ\Re f+h_{3}\circ\Im f)d\mu+i\int_{X}(h_{2}\circ\Re f+h_{4}\circ\Im f)d\mu,\eqno(3.1)$$
then $\Phi$ is a continuous valuation provided that $\delta_{k}=0$ if $\mu(X)=\infty$.
\vskip 5pt
\noindent{\bf Proof.}
For $f,g\in L^{p}(\mathbb{C},\mu)$, let
$$\begin{array}{rl}
\displaystyle &E=\{x\in X:~\Re f\leq \Re g, ~\Im f\leq \Im g\},~~F=\{x\in X:~\Re f\leq \Re g, ~\Im f>\Im g\},
\\\\&G=\{x\in X:~\Re f> \Re g, ~\Im f\leq \Im g\},~~H=\{x\in X:~\Re f> \Re g, ~\Im f> \Im g\}.\end{array}$$
By (3.1) and (1.2), we obtain
$$\begin{array}{rl}
\displaystyle \Phi(f\vee g)&=\displaystyle \int_{X}(h_{1}\circ\Re (f\vee g)+h_{3}\circ\Im (f\vee g))d\mu+i\int_{X}(h_{2}\circ\Re (f\vee g)+h_{4}\circ\Im (f\vee g))d\mu
\\\\&=\displaystyle \int_{E}(h_{1}\circ\Re g+h_{3}\circ\Im g)d\mu+i\int_{E}(h_{2}\circ\Re g+h_{4}\circ\Im g)d\mu
\\\\&~+\displaystyle \int_{F}(h_{1}\circ\Re g+h_{3}\circ\Im f)d\mu+i\int_{F}(h_{2}\circ\Re g+h_{4}\circ\Im f)d\mu
\\\\&~+\displaystyle \int_{G}(h_{1}\circ\Re f+h_{3}\circ\Im g)d\mu+i\int_{G}(h_{2}\circ\Re f+h_{4}\circ\Im g)d\mu
\\\\&~+\displaystyle \int_{H}(h_{1}\circ\Re f+h_{3}\circ\Im f)d\mu+i\int_{H}(h_{2}\circ\Re f+h_{4}\circ\Im f)d\mu.
\end{array}$$
Similarly, by (3.1) and (1.3), we have
$$\begin{array}{rl}
\displaystyle \Phi(f\wedge g)
&=\displaystyle \int_{E}(h_{1}\circ\Re f+h_{3}\circ\Im f)d\mu+i\int_{E}(h_{2}\circ\Re f+h_{4}\circ\Im f)d\mu
\\\\&~+\displaystyle \int_{F}(h_{1}\circ\Re f+h_{3}\circ\Im g)d\mu+i\int_{F}(h_{2}\circ\Re f+h_{4}\circ\Im g)d\mu
\\\\&~+\displaystyle \int_{G}(h_{1}\circ\Re g+h_{3}\circ\Im f)d\mu+i\int_{G}(h_{2}\circ\Re g+h_{4}\circ\Im f)d\mu
\\\\&~+\displaystyle \int_{H}(h_{1}\circ\Re g+h_{3}\circ\Im g)d\mu+i\int_{H}(h_{2}\circ\Re g+h_{4}\circ\Im g)d\mu.
\end{array}$$
 Note that $E\cup F\cup G\cup H=X$ and that $E,F,G,H$ are pairwise disjoint. Thus,
$$
\displaystyle \displaystyle \Phi(f\vee g)+\Phi(f\wedge g)=
\displaystyle \Phi(f)+\Phi(g).$$
Hence $\Phi$ is a valuation.

It remains to show that $\Phi$ is continuous. Let $f\in L^{p}(\mathbb{C},\mu)$ and let $\{f_{k}\}$ be a sequence in $L^{p}(\mathbb{C},\mu)$ with $f_{k}\rightarrow f$ in $L^{p}(\mathbb{C},\mu)$. Next, we will show that $\Phi(f_{k})$ converges to $\Phi(f)$ by showing that every subsequence, $\Phi(f_{k_{l}})$, of $\Phi(f_{k})$ has a subsequence, $\Phi(f_{k_{l_{m}}})$, converges to $\Phi(f)$.
Set $f=\alpha+i\beta$ and $f_{k}=\alpha_{k}+i\beta_{k}$ with $\alpha, \beta, \alpha_{k}, \beta_{k}\in L^{p}(\mathbb{R},\mu)$ such that $\alpha_{k}\rightarrow \alpha$ and $\beta_{k}\rightarrow \beta$ in $L^{p}(\mathbb{R},\mu)$. Let $\{f_{k_{l}}\}$ be a subsequence of $\{f_{k}\}$, then $\{f_{k_{l}}\}$ converges to $f$ in $L^{p}(\mathbb{C},\mu)$. Then there exists a subsequence $\{f_{k_{l_{m}}}\}$ of $\{f_{k_{l}}\}$ which converges to $f$ in $L^{p}(\mathbb{C},\mu)$, where $f_{k_{l_{m}}}=\alpha_{k_{l_{m}}}+i\beta_{k_{l_{m}}}$ with $\alpha_{k_{l_{m}}}, \beta_{k_{l_{m}}}\in L^{p}(\mathbb{R},\mu)$ such that $\alpha_{k_{l_{m}}}\rightarrow \alpha$ and $\beta_{k_{l_{m}}}\rightarrow \beta$ in $L^{p}(\mathbb{R},\mu)$.
Since $h_{1}$ is continuous, we have
$$(h_{1}\circ\alpha_{k_{l_{m}}})(x)\rightarrow (h_{1}\circ\alpha)(x),~a.e.~[\mu].$$
Since $|h_{1}(a)|\leq \gamma_{1}|a|^{p}+\delta_{1}$ for all $a\in\mathbb{R}$, we get
$$|(h_{1}\circ\alpha_{k_{l_{m}}})(x)|\leq \gamma_{1}|\alpha_{k_{l_{m}}}|^{p}+\delta_{1},~a.e.~[\mu].$$
If $\mu(X)<\infty$, apply $\alpha_{k_{l_{m}}}\rightarrow \alpha$ in $L^{p}(\mathbb{R},\mu)$ to get
$$\lim_{m\rightarrow\infty}\int_{X}\gamma_{1}|\alpha_{k_{l_{m}}}|^{p}+\delta_{1}d\mu=\int_{X}\gamma_{1}|\alpha|^{p}d\mu+\delta_{1}\mu(X).$$
And we take $\delta_{1}=0$ in the above equation if $\mu(X)=\infty$.
By a modification of Lebesgue's Dominated Convergence Theorem (see [46, Proposition 2.2.]), we have $h_{1}\circ(\alpha)\in L^{1}(\mathbb{R},\mu)$ and
$$\lim_{m\rightarrow\infty}\int_{X}(h_{1}\circ\alpha_{k_{l_{m}}})d\mu=\int_{X}(h_{1}\circ \alpha) d\mu.\eqno(3.2)$$
Similarly,
$$\lim_{m\rightarrow\infty}\int_{X}(h_{2}\circ\alpha_{k_{l_{m}}})d\mu=\int_{X}(h_{2}\circ \alpha) d\mu,\eqno(3.3)$$
$$\lim_{m\rightarrow\infty}\int_{X}(h_{3}\circ\beta_{k_{l_{m}}})d\mu=\int_{X}(h_{3}\circ \beta) d\mu,\eqno(3.4)$$
$$\lim_{m\rightarrow\infty}\int_{X}(h_{4}\circ\beta_{k_{l_{m}}})d\mu=\int_{X}(h_{4}\circ \beta) d\mu.\eqno(3.5)$$
Thus,
$$\begin{array}{rl}\displaystyle|\Phi(f_{k_{l_{m}}})-\Phi(f)|&=\displaystyle|\int_{X}(h_{1}\circ \alpha_{k_{l_{m}}} -h_{1}\circ \alpha) d\mu+\int_{X}(h_{3}\circ \beta_{k_{l_{m}}} -h_{3}\circ \beta) d\mu
\\\\&~\displaystyle+i\Big(\int_{X}(h_{2}\circ \alpha_{k_{l_{m}}} -h_{2}\circ \alpha) d\mu+\int_{X}(h_{4}\circ \beta_{k_{l_{m}}} -h_{4}\circ \beta) d\mu\Big)|
\\\\&\leq \displaystyle|\int_{X}(h_{1}\circ \alpha_{k_{l_{m}}} -h_{1}\circ \alpha)d\mu|+|\int_{X}(h_{3}\circ \beta_{k_{l_{m}}} -h_{3}\circ \beta) d\mu|
\\\\&~\displaystyle+|\int_{X}(h_{2}\circ \alpha_{k_{l_{m}}} -h_{2}\circ \alpha) d\mu|+|\int_{X}(h_{4}\circ \beta_{k_{l_{m}}} -h_{4}\circ \beta) d\mu|.\end{array}$$
From (3.2)-(3.5), we conclude $\Phi(f_{k_{l_{m}}})\rightarrow \Phi(f)$. Hence, $\Phi$ is continuous. ~~~~~~~~~~$\Box$

\vskip 5pt
\noindent {\bf Lemma 3.2}
If the function $\Phi:L^{p}(\mathbb{C},\mu)\rightarrow \mathbb{C}$ is a valuation, then
$$\Phi(f)=\Phi(\Re f)+\Phi(i \Im f),$$
for all $f\in L^{p}(\mathbb{C},\mu)$.
\vskip 5pt
\noindent{\bf Proof.}
If $\Re f, \Im f\geq 0$ or $\Re~f, \Im~f\leq 0$, then, by (1.2) and (1.3), we have
$$\Phi(\Re f)+\Phi(i \Im f)=\Phi(f)+\Phi(0).\eqno(3.6)$$
If $\Re f\geq 0, \Im f\leq 0$ or $\Re f\leq 0, \Im f\geq 0$, then, by (1.2) and (1.3), we have
$$\Phi(f)+\Phi(0)=\Phi(\Re f)+\Phi(i \Im f).\eqno(3.7)$$
Note that $\Phi(0)=0$, apply (3.6) and (3.7), get
$$\Phi(f)=\Phi(\Re f)+\Phi(i \Im f)$$
for all $f\in L^{p}(\mathbb{C},\mu)$.~~~~~~~~~~~~~~~~~~~~~~~~~~~~~~~~~~~~~~~~~~~~~~~~~~~~~~~~~~~~~~~~~~~~~~~~~~~~~~~~~~~~~~$\Box$
\vskip 5pt
If we restrict $f$ to belong to $L^{p}(\mathbb{R},\mu)$, then it is obvious that $\Phi$ is a valuation on $L^{p}(\mathbb{R},\mu)$.
Also, we can construct another valuation on $L^{p}(\mathbb{R},\mu)$ which is related to $\Phi$.
\vskip 5pt
\noindent {\bf Lemma 3.3}
Let $\Phi: L^{p}(\mathbb{C},\mu)\rightarrow \mathbb{C}$ be a valuation. If the functions $\Phi': L^{p}(\mathbb{R},\mu)\rightarrow \mathbb{C}$ is defined by
$$\Phi'(f)=\Phi (if)$$
for all $f\in L^{p}(\mathbb{R},\mu)$, then
$\Phi'$ is a valuation on $L^{p}(\mathbb{R},\mu)$.

\vskip 5pt
\noindent{\bf Proof.}
For $f,g\in L^{p}(\mathbb{R},\mu)$, by (1.2) and (1.3), we have
$$i(f\vee g)=if\vee ig~~\textrm{and}~~i(f\wedge g)=if\wedge ig.\eqno(3.8)$$
By (3.8) and the valuation property of $\Phi$, it follows that
$$\Phi'(f\vee g)+\Phi'(f\wedge g)=\Phi(if\vee ig)+\Phi(if\wedge ig)=\Phi(if)+\Phi(ig)=\Phi'(f)+\Phi'(g)$$
for all $f,g\in L^{p}(\mathbb{R},\mu)$. Thus, $\Phi'$ is a valuation on $L^{p}(\mathbb{R},\mu)$. ~~~~~~~~~~~~~~~~~~~~~~~~~~$\Box$

\vskip 5pt
\noindent {\bf Lemma 3.4}
Let $\Phi: L^{p}(\mathbb{R},\mu)\rightarrow \mathbb{C}$ be a valuation. If the functions $\Phi_{1}, \Phi_{2}: L^{p}(\mathbb{R},\mu)\rightarrow \mathbb{R}$ are defined by
$$\Phi(f)=\Phi_{1}(f)+i\Phi_{2}(f)$$
for all $f\in L^{p}(\mathbb{R},\mu)$, then
$\Phi_{1}, \Phi_{2}$ both are a real-valued valuation on $L^{p}(\mathbb{R},\mu)$.

\vskip 5pt
\noindent{\bf Proof.} Since $\Phi$ is a valuation, we have
$$\begin{array}{rl}\displaystyle &\Phi(f\vee g)+\Phi(f\wedge g)=\Phi_{1}(f\vee g)+i\Phi_{2}(f\vee g)+\Phi_{1}(f\wedge g)+i\Phi_{2}(f\wedge g)
\\\\&=\displaystyle\Phi(f)+\Phi(g)=\Phi_{1}(f)+i\Phi_{2}(f)+\Phi_{1}(g)+i\Phi_{2}(g)\end{array}$$
for all $f,g\in L^{p}(\mathbb{R},\mu)$. Thus,
$$\Phi_{1}(f\vee g)+\Phi_{1}(f\wedge g)=\Phi_{1}(f)+\Phi_{1}(g),$$
and
$$\Phi_{2}(f\vee g)+\Phi_{2}(f\wedge g)=\Phi_{2}(f)+\Phi_{2}(g),$$
for all $f,g\in L^{p}(\mathbb{R},\mu)$. Therefore, $\Phi_{1}, \Phi_{2}$ both are a real-valued valuation on $L^{p}(\mathbb{R},\mu)$.
~~~~~~~~~~~~~~~~~~~~~~~~~~~~~~~~~~~~~~~~~~~~~~~~~~~~~~~~~~~~~~~~~~~~~~~~~~~~~~~~~~~~~~~~~~~$\Box$

\vskip 5pt
In order to establish a representation theorem for continuous complex-valued valuations on $L^{p}(\mathbb{C},\mu)$, we will use the corresponding representation theorem for real case which was obtained by Tsang [46].
\vskip 5pt \noindent {\bf Theorem 3.5} ([46]) Let $(X,\mathfrak{F},\mu)$ be a non-atomic measure space and let $\Phi: L^{p}(\mathbb{R},\mu)\rightarrow \mathbb{R}$ be a continuous translation invariant valuation. If
there exists a continuous function $h:\mathbb{R}\rightarrow \mathbb{R}$ with $h(0)=0$ such that
$\Phi(b\chi_{E})=h(b)\mu(E)$
for all $b\in \mathbb{R}$ and all $E\in\mathfrak{F}$ with $\mu(E)<\infty$, then there exist constants $\gamma,\delta\geq 0$
such that $|h(a)|\leq \gamma|a|^{p}+\delta$ for all $a\in \mathbb{R}$, and
$$\Phi(f)=\int_{X}(h\circ f)d\mu$$
for all $f\in L^{p}(\mathbb{R},\mu)$. In addition, $\delta=0$ if $\mu(X)=\infty$.

\vskip 5pt
\noindent {\bf Proof of Theorem 1.1}
\vskip 5pt Let $\Phi: L^{p}(\mathbb{C},\mu)\rightarrow \mathbb{C}$ be a continuous valuation. For $f\in L^{p}(\mathbb{C},\mu)$,
by Lemma 3.2, Lemma 3.3 and Lemma 3.4, we have
$$\begin{array}{rl}\displaystyle\Phi(f)&=\Phi(\Re f)+\Phi(i\Im f)=\Phi_{1}(\Re f)+i\Phi_{2}(\Re f)+\Phi_{1}(i\Im f)+i\Phi_{2}(i\Im f)
\\\\&=\displaystyle\Phi_{1}(\Re f)+i\Phi_{2}(\Re f)+\Phi_{1}'(\Im f)+i\Phi_{2}'(\Im f),\end{array}$$
where $\Phi(\Re f)=\Phi_{1}(\Re f)+i\Phi_{2}(\Re f), \Phi_{1}'(\Im f)=\Phi_{1}(i\Im f), ~\textrm{and}~\Phi_{2}'(\Im f)=i\Phi_{2}(i\Im f)$.
Since $\Re (f\vee g)=\Re f\vee \Re g, \Re (f\wedge g)=\Re f\wedge \Re g$, $\Im (f\vee g)=\Im f\vee \Im g$, and $\Im (f\wedge g)=\Im f\wedge \Im g$.
Moreover, Lemma 3.3 and Lemma 3.4 imply that $\Phi_{1},\Phi_{2},\Phi_{1}', \Phi_{2}'$ all are a real-valued valuation on $L^{p}(\mathbb{R},\mu)$.

\vskip 5pt
If we restrict to $f\in L^{p}(\mathbb{R},\mu)$, then the continuity of $\Phi$ implies that $\Phi_{1},\Phi_{2}$ are continuous on $L^{p}(\mathbb{R},\mu)$.
If we consider $f\in L^{p}(\mathbb{C},\mu)$ with $\Re f=0$, then the continuity of $\Phi$ implies that $\Phi_{1}',\Phi_{2}'$ are continuous on $L^{p}(\mathbb{R},\mu)$. Thus, $\Phi_{1},\Phi_{2},\Phi_{1}', \Phi_{2}'$ all are a continuous real-valued valuation on $L^{p}(\mathbb{R},\mu)$. It follows from Theorem 3.5 that there
exist continuous functions $h_{k}:\mathbb{R}\rightarrow \mathbb{R}$ with the properties that $h_{k}(0)=0$ and there exist constants $\gamma_{k},\delta_{k}\geq 0$ such that
$|h_{k}(a)|\leq \gamma_{k}|a|^{p}+\delta_{k}$ for all $a\in\mathbb{R}~(k=1,2,3,4)$, and
$$\Phi(f)=\int_{X}(h_{1}\circ\Re f+h_{3}\circ\Im f)d\mu+i\int_{X}(h_{2}\circ\Re f+h_{4}\circ\Im f)d\mu$$
for all $f\in L^{p}(\mathbb{C},\mu)$. In addition, $\delta_{k}=0~(k=1,2,3,4)$ if $\mu(X)=\infty$.
~~~~~~~~~~~~~~~~~~~~~~~~$\Box$

\vskip 5pt
If $\mu$ is Lebesgue measure, then $L^{p}(\mathbb{C},\mu)$ becomes the space of Lebesgue measurable complex-valued functions.
We usually write as $L^{p}(\mathbb{C},\mathbb{R}^{n})$.
\vskip 5pt
\noindent {\bf Lemma 3.6} Let $h_{k}: \mathbb{R}\rightarrow \mathbb{R}$ be continuous functions with the property that
there exist $\gamma_{k}\geq 0$ such that $|h_{k}(a)|\leq \gamma_{k}|a|^{p}$ for all $a\in\mathbb{R}~(k=1,\ldots,4)$. If the function $\Phi: L^{p}(\mathbb{C},\mathbb{R}^{n})\rightarrow \mathbb{C}$ is defined by
$$\Phi(f)=\int_{\mathbb{R}^{n}}(h_{1}\circ\Re f+h_{3}\circ\Im f)dx+i\int_{X}(h_{2}\circ\Re f+h_{4}\circ\Im f)dx,$$
then $\Phi$ is a continuous translation invariant valuation.

\vskip 5pt
\noindent{\bf Proof.} Let $\mathcal{M}$ denote the collection of Lebesgue measurable sets in $\mathbb{R}^{n}$. Take $X=\mathbb{R}^{n}$, $\mathfrak{F}=\mathcal{M}$ and $\mu$ Lebesgue measure in Lemma 3.1 to conclude that $\Phi$ is a continuous valuation on $L^{p}(\mathbb{C},\mathbb{R}^{n})$.

For every $t\in \mathbb{R}^{n}$ and every $f\in L^{p}(\mathbb{C},\mathbb{R}^{n})$, we have
$$\displaystyle \Phi(f(x-t))=\displaystyle \int_{\mathbb{R}^{n}}(h_{1}\circ\Re f+h_{3}\circ\Im f)(x-t)dx+i\int_{X}(h_{2}\circ\Re f+h_{4}\circ\Im f)(x-t)dx=\displaystyle \Phi(f),$$
which means that $\Phi$ is translation invariant. ~~~~~~~~~~~~~~~~~~~~~~~~~~~~~~~~~~~~~~~~~~~~~~~$\Box$

\vskip 15pt
\noindent{\bf Proof of Theorem 1.2}
 \vskip 5pt
It follows from Lemma 3.6 that (1.4) determines a continuous translation invariant valuation on $L^{p}(\mathbb{C},\mathbb{R}^{n})$.

Conversely, let $\Phi: L^{p}(\mathbb{C},\mathbb{R}^{n})\rightarrow \mathbb{C}$ be a continuous translation invariant valuation.
Taking $X=\mathbb{R}^{n}$, $\mathfrak{F}=\mathcal{M}$
and $\mu$ Lebesgue measure in the proof of Theorem 1.1, we obtain
$$\Phi(f)=\displaystyle\Phi_{1}(\Re f)+i\Phi_{2}(\Re f)+\Phi_{1}'(\Im f)+i\Phi_{2}'(\Im f),$$
where $\Phi_{1},\Phi_{2},\Phi_{1}', \Phi_{2}'$ all are a real-valued valuation on $L^{p}(\mathbb{R},\mu)$. Theorem 2.1 implies that
there exist continuous functions $h_{k}:\mathbb{R}\rightarrow \mathbb{R}$ with the property that there exist constants $\gamma_{k}\geq 0$ such that
$|h_{k}(a)|\leq \gamma_{k}|a|^{p}$ for all $a\in\mathbb{R}~(k=1,2,3,4)$, and
$$\Phi(f)=\int_{\mathbb{R}^{n}}(h_{1}\circ\Re f+h_{3}\circ\Im f)dx+i\int_{\mathbb{R}^{n}}(h_{2}\circ\Re f+h_{4}\circ\Im f)dx$$
for all $f\in L^{p}(\mathbb{C},\mathbb{R}^{n})$.
~~~~~~~~~~~~~~~~~~~~~~~~~~~~~~~~~~~~~~~~~~~~~~~~~~~~~~~~~~~~~~~~~~~~~~~~~~~~~~$\Box$

 \vskip 5pt

Let $\mathcal{W}$ denote the $\sigma$-algebra defined as
$$\mathcal{W}=\{E: E\subseteq S^{n-1}, \{\lambda x: x\in E, 0\leq \lambda\leq 1\}\in\mathcal{M}\}.$$
Also denote by $\sigma$ the spherical Lebesgue measure.
If $\mu$ is the spherical Lebesgue measure, then $L^{p}(\mathbb{C},\sigma)$ denotes the space of spherical Lebesgue measurable complex-valued functions. We usually write as $L^{p}(\mathbb{C},S^{n-1})$.
\vskip 5pt
\noindent {\bf Lemma 3.7} Let $h_{k}: \mathbb{R}\rightarrow \mathbb{R}$ be continuous functions with the properties that $h_{k}(0)=0$ and
there exist $\gamma_{k}, \delta_{k}\geq 0$ such that $|h_{k}(a)|\leq \gamma_{k}|a|^{p}+\delta_{k}$ for all $a\in\mathbb{R}~(k=1,\ldots,4)$. If the function $\Phi: L^{p}(\mathbb{C},S^{n-1})\rightarrow \mathbb{C}$ is defined by
$$\Phi(f)=\int_{S^{n-1}}(h_{1}\circ\Re f+h_{3}\circ\Im f)du+i\int_{S^{n-1}}(h_{2}\circ\Re f+h_{4}\circ\Im f)du,$$
then $\Phi$ is a continuous rotation invariant valuation.

\vskip 5pt
\noindent{\bf Proof.}
 Take $X=S^{n-1}$, $\mathfrak{F}=\mathcal{W}$
and $\mu=\sigma$ in Lemma 3.1 to conclude that $\Phi$ is a continuous valuation on $L^{p}(\mathbb{C},S^{n-1})$.

Note that $\theta u\in S^{n-1}$ for every $\theta\in \textrm{O}(n)$ and every $u\in S^{n-1}$. Since the spherical Lebesgue measure is rotation invariant, we have
$$\displaystyle \Phi(f\circ \theta^{-1})=\displaystyle \int_{S^{n-1}}(h_{1}\circ\Re f+h_{3}\circ\Im f)(\theta u)du+i\int_{S^{n-1}}(h_{2}\circ\Re f+h_{4}\circ\Im f)(\theta u)du=\displaystyle \Phi(f)$$
for all $f\in L^{p}(\mathbb{C},S^{n-1})$,
which finishes the proof.~~~~~~~~~~~~~~~~~~~~~~~~~~~~~~~~~~~~~~~~~~~~~~~~~~~~~~~~~~~~~~~~~~$\Box$

\vskip 15pt
\noindent{\bf Proof of Theorem 1.3}
\vskip 5pt
 \vskip 5pt
It follows from Lemma 3.7 that (1.5) determines a continuous rotation invariant valuation on $L^{p}(\mathbb{C},S^{n-1})$.

Conversely, let $\Phi: L^{p}(\mathbb{C},S^{n-1})\rightarrow \mathbb{C}$ be a continuous rotation invariant valuation.
Taking $X=S^{n-1}$, $\mathfrak{F}=\mathcal{W}$ and $\mu=\sigma$ in the proof of Theorem 1.1, we obtain
$$\Phi(f)=\displaystyle\Phi_{1}(\Re f)+i\Phi_{2}(\Re f)+\Phi_{1}'(\Im f)+i\Phi_{2}'(\Im f),$$
where $\Phi_{1},\Phi_{2},\Phi_{1}', \Phi_{2}'$ all are a real-valued valuation on $L^{p}(\mathbb{R},S^{n-1})$. Theorem 2.2 implies that
there exist continuous functions $h_{k}:\mathbb{R}\rightarrow \mathbb{R}$ with the properties that $h_{k}(0)=0$ and there exist constants $\gamma_{k}, \delta_{k}\geq 0$ such that
$|h_{k}(a)|\leq \gamma_{k}|a|^{p}+\delta_{k}$ for all $a\in\mathbb{R}~(k=1,2,3,4)$, and
$$\Phi(f)=\int_{S^{n-1}}(h_{1}\circ\Re f+h_{3}\circ\Im f)du+i\int_{S^{n-1}}(h_{2}\circ\Re f+h_{4}\circ\Im f)du$$
for all $f\in L^{p}(\mathbb{C},S^{n-1})$.
~~~~~~~~~~~~~~~~~~~~~~~~~~~~~~~~~~~~~~~~~~~~~~~~~~~~~~~~~~~~~~~~~~~~~~~~~~~~~$\Box$

\vskip 10pt
\noindent {\bf Acknowledgment} The work was supported in part by the Natural Science Foundation of Hunan Province (2019JJ50172).

\vskip 30pt
\end{document}